\documentclass[a4paper,11pt]{article}

\usepackage{mathrsfs}

\usepackage{amsfonts}
\usepackage{graphicx}

\usepackage{amssymb}
\usepackage{amssymb,amsmath,amsthm}

\usepackage{graphicx}
\usepackage{amsmath}
\usepackage{amssymb,amsmath,amsthm}
\usepackage{indentfirst}

\numberwithin{equation} {section} \makeatletter
\renewcommand{\@seccntformat}[1]{\csname the#1\endcsname.\hspace{0.5em}} \makeatother

\begin{document}

\title{\textbf{Concentration in  flux-function limits  of solutions to a deposition model}
 \thanks{\ Supported by Applied Basic Research Projects of Yunnan Province (2015FB104).}}

\author{Hongjun Cheng\ \ \ \ \ Shiwei Li\\
 {\small{ Department of Mathematics, Yunnan University,  Kunming, Yunnan 650091, People's Republic of China }}
       }

\date{}
\maketitle

\noindent\textbf{Abstract}\ This paper is concerned with a singular
flux-function limit of the Riemann solutions to a deposition model.
As a result, it is shown that the Riemann solutions to the
deposition model just converge to the corresponding Riemann
solutions to the limit system, which is one of typical models
admitting delta-shocks. Especially, the phenomenon of concentration
and the formation of delta-shocks in the limit are  analyzed in
detail, and the process  of concentration is numerically simulated.

\vspace{0.2cm}

\noindent\textbf{keywords}\ Deposition model; Flux-function limits;
Concentration; Delta-shocks; Numerical simulations.

\vspace{0.2cm}


\section{Introduction}

We consider the following deposition model
\begin{equation}
\left\{\begin{array}{l}
 v_t+(uv)_x=0,\\
 u_t+ (u^2+\epsilon v)_x=0,
\end{array}
\right.\label{DM}
\end{equation}
where $ v\geq0$ is the density of the population performing the
deposition, $u=-\partial_x h$ with $h=h(x,t)$ being the deposition
height, and $\epsilon$ is a positive parameter. The first equation
describes the conservation of total population. The second one is
derived from the rules governing the time evolution of the
deposition system: (1) the deposition consists of
population-generating deposition and self-generating deposition; (2)
the population is driven by a velocity field proportional to the
negative gradient of height. Besides, model (\ref{DM}) can also be
derived as  the hydrodynamic limit of  some randomly growing
interface models \cite{F-T, F-V}. The system  (\ref{DM})  is also
called as the Leroux system in the PDE literature \cite{Le, Se-1}.
For some investigations concerning (\ref{DM}), see
\cite{Se-2,Gr,L-M,Lu}, etc.

One can notice that as $\epsilon\rightarrow0^+$,  the model
(\ref{DM}) formally becomes  the following system
\begin{equation}
\left\{\begin{array}{l}
 u_t+ (u^2)_x=0,\\
  v_t+( v u)_x=0.
\end{array}\right.\label{Bur}
\end{equation}
This is one of very typical models in the literature with respect to
delta-shocks
\cite{B,SZh,Tan-Zhang,Fl,KK,K,Yang-1,DS,VSh,Cheng-Yang-1,Cheng-1},
an interesting topic. It is a mathematical simplification of Euler
equations of gas dynamics and can be obtained by setting density and
pressure to be constant in the momentum conservation laws. It also
has some physical interpretations. For instance, it can  be used to
model the flow of particles with $u$ being the velocity and $v$ the
density. In 1977, Korchinski \cite{Ko} considered (\ref{Bur}) in his
unpublished Ph.D. thesis. Motivated by his numerical study, he used
a kind of generalized delta-function in the construction of his
unique solution to the Riemann problem. Afterwards, in 1994, Tan et
al. \cite{T-Z-Z}  found that in the Riemann problem for (\ref{Bur}),
no classical weak solution exists and delta-shocks should be
introduced for some initial data. With the delta-shocks, they solved
the Riemann problem completely.  Under some reasonable second order
viscous approximations, the stability of delta-shocks for
(\ref{Bur}) was also proved in \cite{T-Z-Z, Ta, Hu}.

The main purpose of this paper is to study the  behaviors of
solutions of system (\ref{DM}) as the flux  $\epsilon v$ vanishes
(that is, $\epsilon\rightarrow0^+$) by the  Riemann problem. We are
especially concerned with the phenomenon of concentration and the
formation of delta-shocks in the limit.

Firstly, we consider the Riemann problem  for  (\ref{DM}) with
initial data
\begin{equation}
(u,v)(x,t=0)=(u_\pm,v_\pm),\ \ \ \ \pm x>0,\label{I-D}
\end{equation}
where $v_\pm>0$. System (\ref{DM}) is nonstrictly hyperbolic, and
both characteristic fields are genuinely nonlinear.  The elementary
waves include shocks and rarefaction waves, and  (\ref{DM}) belongs
to the so-called Temple class \cite{Te}. By the analysis method in
phase plane,  the unique global Riemann solution is constructed with
four different kinds of structures containing  shock(s) and/or
rarefaction wave(s).

Secondly, we study the  behaviors of solutions of system (\ref{DM})
as the flux  $\epsilon v$ vanishes by the  Riemann problem. As a
result, it is rigorously shown  that as $\epsilon\rightarrow0^+$,
the Riemann solutions to (\ref{DM}) just converge to the Riemann
solutions to (\ref{Bur}) with the same initial data. Especially,
when $u_ +\leq0\leq u_-$, the two-shock solution to (\ref{DM}) and
(\ref{I-D}) tends to the delta-shock solution to (\ref{Bur}) and
(\ref{I-D}), where the intermediate density between the two shocks
tends to a weighted $\delta$-measure which forms the delta-shock.
Further,  the process of concentration is numerically simulated. It
can be seen that such a  flux-function limit may be very singular:
the limit functions of solutions are no longer in the spaces of
functions $BV$ or $L^\infty$, and  the space of Radon measures, for
which the divergences of certain entropy and entropy flux fields are
also Radon measures, is a natural space in order to deal with such a
limit.

Let us remark that in the past more than 10 years, more attention
has been paid on the investigation of phenomenon of concentration
and the formation of delta-shocks in solutions to hyperbolic systems
of conservation laws. Li \cite{LJQ} and  Chen and Liu \cite{C-L-1,
C-L-2} identified and analyzed  the phenomenon of concentration and
the formation of delta-shocks in solutions to the Euler equations
for both isentropic and nonisentropic fluids as the pressure
vanishes. Yin and Sheng \cite{YSh-1,YSh-2} extended the studies to
the relativistic Euler equations.   With respect to this topic, also
see \cite{D-MN,Yang-Liu-1,Yang-Liu-2,Sh}.

We arrange the rest of the paper as follows. In the following
section, we recall the Riemann problem for system (\ref{Bur}). In
Section 3, we solve the Riemann problem for  (\ref{DM})  by the
analysis method in phase-plane. In Section 4 and Section 5, we study
the limits  of solutions of the Riemann problem for (\ref{DM})  as
$\epsilon\rightarrow0^+$. In   Section 6,  we  examine the process
of concentration as $\epsilon$ decreases by some numerical results.

\section{Solutions of the Riemann problem for (\ref{Bur})}

In this section, we  recall  the Riemann problem for (\ref{Bur})
with initial data (\ref{I-D})  which was solved by Tan et al.
\cite{T-Z-Z}.  The characteristic roots of (\ref{Bur}) are
$\lambda_1=u$ and $\lambda_2=2u$,  and the corresponding right
characteristic vectors are $r_1=(0,1)^T$ and $r_2=(1,v/u)^T$,
respectively. They satisfy $\nabla\lambda_1\cdot r_1\equiv0$ and
$\nabla\lambda_2\cdot r_2=2$, where and in the following
$\nabla=({\partial}/{\partial u},{\partial}/{\partial v})$ is the
gradient operator. Therefore (\ref{Bur}) is nonstrictly hyperbolic
because of $\lambda_1=\lambda_2$ at $u=0$, $\lambda_1$ is linearly
degenerate, and $\lambda_2$ is genuinely nonlinear.

 Since the equations and the Riemann data are invariant under
uniform stretching of coordinates $(x,t)\rightarrow(\beta x,\beta
t)(\beta>0)$, we consider the self-similar solutions  $(u,
v)(x,t)=(u, v)(\xi)$, where $\xi=x/t$. Then the Riemann problem
turns into
\begin{equation}
\left\{\begin{array}{l}
 -\xi u_\xi+(u^2)_\xi=0,\cr\noalign{\vskip2truemm}
 -\xi v_\xi+( v u)_\xi=0,
 \end{array}\right.\label{2.1}
\end{equation}
and
\begin{equation}
(u, v)(\pm\infty)=(u_\pm, v_\pm).\label{2.2}
\end{equation}
This is a two-point boundary value problem of first-order ordinary
differential equations with the boundary values in the infinity.

Besides the constant states, the self-similar waves
$(u,v)(\xi)(\xi=x/t)$ of the first family are contact
discontinuities
\begin{equation}
J:\
 \xi=u_l=u_r,\label{2.3}
\end{equation}
and those of the second family are rarefaction waves
\begin{equation}
R:\
 \xi=2u,\ \ u/v=u_l/v_l,\ \ u>u_l,
\label{2.4}
\end{equation}
or shocks
\begin{equation}
S:\
 \xi=u_l+u_r,\ \
 u_r/v_r=u_l/v_l,\  \ u_l>u_r>0\  \  \mbox{or}\ \  0>u_l>u_r,
 \label{2.5}
\end{equation}
where the indices $l$ and $r$ denote the left and right states
respectively. All of $J, R$ and $S$ are waves with
$(u(\xi),v(\xi))\in BV$ and are called the classical waves.

Using the classical waves, by the analysis in phase-plane, one can
construct the solutions of Riemann problem (\ref{Bur}) and
(\ref{I-D}) in the following cases
$$
\begin{array}{lll}
 R+J\ (u_-<u_+<0),&J+R\ (0<u_-<u_+),&R+R\ (u_-\leq0\leq u_+),\\[2mm]
  S+J\ (u_+<u_-<0),& J+S\ (0<u_+<u_-).&
\end{array}
$$
However, for the case $u_+\leq0\leq u_-$,  the singularity cannot be
a jump with finite amplitude; that is, there is no solution which is
piecewise smooth and bounded. Hence a solution containing a weighted
$\delta$-measure (i.e., delta-shock) supported on a line should be
introduced  in order to establish the existence in a space of
measures from the mathematical point of view.

We define the weighted  $\delta$-measure $w(s)\delta_L$ supported on
a smooth curve $L$ parameterized as $t=t(s), x=x(s)\ (c\leq s\leq
d)$  by
\begin{equation}
\Big<w(s)\delta_L,\psi(x,t)\Big>=\displaystyle\int_c^dw(s)\psi\big(t(s),x(s)\big)ds\label{2.6}
\end{equation}
for all test functions $\psi(x,t)\in
C^\infty_0((-\infty,+\infty)\times[0,+\infty))$.

With this definition,  when $u_+\leq0\leq u_-$, the solution of
Riemann problem  (\ref{Bur})  and (\ref{I-D}) is the following
solution involving delta-shock  in the form
\begin{equation}
(u,v)(x,t)=\left\{\begin{array}{ll}
  (u_-,v_-),&x<x(t),\\[1mm]
  \big(u_\delta(t), w(t)\delta(x-x(t))\big),&x=x(t),\\[1mm]
  (u_+,v_+),&x>x(t)
\end{array}\right.\label{2.7}
\end{equation}
satisfying the  generalized Rankine-Hugoniot relation
\begin{equation}
\left\{\begin{array}{lll}
 \displaystyle \frac{dx(t)}{dt}=u_\delta(t)\\[2mm]
  -u_\delta(t)[u]+[u^2]=0,\\[2mm]
  \displaystyle\frac{dw(t)}{dt}=-u_\delta(t)[v]+[uv]
\end{array}\right.\label{2.8}
\end{equation}
and the entropy condition
\begin{equation}
\lambda_2(u_+)\leq\lambda_1(u_+)\leq
u_\delta(t)\leq\lambda_1(u_-)\leq\lambda_2(u_-),\label{2.9}
\end{equation}
where $[g]=g_--g_+$  is the jump of $g$ across the discontinuity.
Solving the generalized Rankine-Hugoniot relation (\ref{2.8}) under
the entropy condition (\ref{2.9}) gives
\begin{equation}
\left\{\begin{array}{lll}
 x(t)=(u_-+u_+)t,\\
 u_\delta(t)=u_-+u_+,\\
 w(t)=(u_-v_+-u_+v_-)t.
\end{array}\right.\label{2.10}
\end{equation}

\section{Solutions of the Riemann problem for (\ref{DM})}

In this section, we solve the Riemann problem for system (\ref{DM})
with initial data $(\ref{I-D})$, and examine the dependence of the
Riemann solutions on the parameter $\epsilon>0$. Also see the paper
\cite{J-S}. The characteristic roots and corresponding right
characteristic vectors of (\ref{DM}) are
$$
\begin{array}{l}
\lambda_1^\epsilon=u+\displaystyle\frac{u-\sqrt{u^2+4\epsilon
v}}{2},\ \ \ \
\lambda_2^\epsilon=u+\displaystyle\frac{u+\sqrt{u^2+4\epsilon
v}}{2},
\end{array}
$$

$$
\begin{array}{l}
\overrightarrow{r}_1^\epsilon=\bigg(1,\displaystyle\frac{-u-\sqrt{u^2+4\epsilon
v}}{2\epsilon }\bigg)^T, \ \ \ \
\overrightarrow{r}_2^\epsilon=\bigg(1,\displaystyle\frac{-u+\sqrt{u^2+4\epsilon
v}}{2\epsilon }\bigg)^T.
\end{array}
$$
It is easy to calculate  $\nabla\lambda_i^\epsilon\cdot
\overrightarrow{r}_i^\epsilon=2\ (i=1,2)$. So (\ref{DM}) is
nonstrictly hyperbolic, both characteristic fields are genuinely
nonlinear. Moreover, the Riemann invariants along with the
characteristic fields may be selected as, respectively,
\begin{equation}
w(u,v)=\displaystyle\frac{-u-\sqrt{u^2+4\epsilon v}}{2\epsilon},\ \
\ \ z(u,v)=\displaystyle\frac{-u+\sqrt{u^2+4\epsilon
v}}{2\epsilon}.\label{3.1}
\end{equation}

As usual, we seek the self-similar solutions $(u, v)(x,t)=(u,
v)(\xi)$, where $\xi=x/t$. Then the Riemann problem becomes the
boundary value problem
\begin{equation}
\left\{\begin{array}{l}
 -\xi u_\xi+\big(u^2+\epsilon v\big)_\xi=0,\\
 -\xi v_\xi+( v u)_\xi=0,
\end{array}\right.\label{3.2}
\end{equation}
and
\begin{equation}
(u, v)(\pm\infty)=(u_\pm, v_\pm).\label{3.3}
\end{equation}

For any smooth solution, (\ref{3.2}) becomes
\begin{equation}
 \left(\begin{array}{cc}
 2u-\xi&\epsilon \\
 v&u-\xi
\end{array}\right)
\left(\begin{array}{l}
 u\\
 v
\end{array}\right)_\xi
=0. \label{3.4}
\end{equation}
Besides  the constant states, the smooth solutions are composed of
the 1-rarefaction waves
\begin{equation}
\left\{\begin{array}{l}
 \xi=\lambda_1^\epsilon=u+\displaystyle\frac{u-\sqrt{u^2+4\epsilon v}}{2},\cr\noalign {\vskip4truemm}
  u-u_0=\Bigg(\displaystyle\frac{u_0-\sqrt{u_0^2+4\epsilon v_0}}{2 v_0}\ \Bigg)( v- v_0),
\end{array}\right.\label{3.5}
\end{equation}
and the 2-rarefaction waves
\begin{equation}
\left\{\begin{array}{l}
 \xi=\lambda_2^\epsilon=u+\displaystyle\frac{u+\sqrt{u^2+4\epsilon v}}{2},\cr\noalign {\vskip2truemm}
   u- u_0=\Bigg(\displaystyle\frac{u_0+\sqrt{u_0^2+4\epsilon
   v_0}}{2 v_0}\Bigg)(v-v_0),
\end{array}\right.\label{3.6}
\end{equation}
where $(u_0,v_0)$ is any state. For them, we have
\begin{equation}
\displaystyle\frac{\mbox{d}\lambda_i^\epsilon}{\mbox{d}u}=\displaystyle\frac{\partial\lambda_i^\epsilon}{\partial
u}+\displaystyle\frac{\partial\lambda_i^\epsilon}{\partial
v}\frac{\mbox{d} v}{\mbox{d}u}=2>0,\ \ \ i=1,2.\label{3.7}
\end{equation}

Let $(u_l,v_l)$ and $(u_r,v_r)$ denote the  states connected by a
rarefaction wave on the left and right sides respectively. Then the
condition $\lambda_1^\epsilon(u_r,v_r)>\lambda_1^\epsilon(u_l,v_l)$
and $\lambda_2^\epsilon(u_r,v_r)>\lambda_2^\epsilon(u_l,v_l)$  are
required for the 1- and 2-rarefaction wave, respectively. From
(\ref{3.7}), it is known that both the 1- and 2-rarefaction wave
should satisfy
\begin{equation}
u_r>u_l.\label{3.8}
\end{equation}

For a given state $(u_l,v_l)$, all possible states which can connect
to $(u_l,v_l)$ on the right by a 1-rarefaction wave must be located
on the straight line
\begin{equation}
{R_1(u_l,v_l)}:\ \ \ \ \
u-u_l=\Bigg(\displaystyle\frac{u_l-\sqrt{u_l^2+4\epsilon v_l}}{2
v_l}\ \Bigg)( v- v_l),\ \
 \ \ u>u_l,
\label{3.9}
\end{equation}
and all possible states which can  connect  to $(u_l,v_l)$ on the
right by a 2-rarefaction wave must be located on straight line
\begin{equation}
{R_2(u_l,v_l)}:\ \ \ \ \
 u- u_l=\Bigg(\displaystyle\frac{u_l+\sqrt{u_l^2+4\epsilon
   v_l}}{2 v_l}\Bigg)(v-v_l),\ \
 \ \ u>u_l.
\label{3.10}
\end{equation}

Let us turn to the discontinuous solutions. For a bounded
discontinuity  at  $x=x(t)$,  the Rankine-Hugoniot  relation reads
\begin{equation}
\left\{\begin{array}{l}
 -\sigma[u]+[ u^2+\epsilon v]=0,\cr\noalign {\vskip1truemm}
 -\sigma[ v]+[ v u]=0,
\end{array}\right.\label{3.11}
\end{equation}
where $\sigma=dx/dt$, $[u]=u_l-u_r$ with $u_l=u(x(t)-0,t)$ and
$u_r=u(x(t)+0,t)$, and so forth.

From (\ref{3.11}), one easily obtains
\begin{equation}
\epsilon\bigg(\frac{[v]}{[u]}\bigg)^2+(u_l+u_r)\frac{[v]}{[u]}-\frac{[uv]}{[u]}=0.\label{3.12}
\end{equation}
By noticing
$$
\frac{[uv]}{[u]}=v_l+u_r\frac{[v]}{[u]},
$$
we solve (\ref{3.12}) to obtain
\begin{equation}
\frac{[v]}{[u]}=\frac{-u_l\pm\sqrt{u_l^2+4\epsilon
v_l}}{2\epsilon}.\label{3.13}
\end{equation}
Then we obtain two kinds of discontinuities
\begin{equation}
\left\{\begin{array}{l}
 \sigma_1=u_r+\displaystyle\frac{u_l-\sqrt{u_l^2+4\epsilon v_l}}{2},\cr\noalign {\vskip2truemm}
 u_r-u_l=\Bigg(\displaystyle\frac{u_l-\sqrt{u_l^2+4\epsilon v_l}}{2 v_l}\Bigg)(v_r- v_l)
\end{array}\right.\label{3.14}
\end{equation}
and
\begin{equation}
 \left\{\begin{array}{l}
 \sigma_2=u_r+\displaystyle\frac{u_l+\sqrt{u_l^2+4\epsilon v_l}}{2},\cr\noalign{\vskip2truemm}
  u_r- u_l=\Bigg(\displaystyle\frac{u_l+\sqrt{u_l^2+4\epsilon v_l}}{2
  v_l}\Bigg)(v_r-v_l).
\end{array}\right.\label{3.15}
\end{equation}
Notice that the second equations in (\ref{3.14}) and (\ref{3.15})
are equivalent to
\begin{equation}
\frac{-u_l-\sqrt{u_l^2+4\epsilon
v_l}}{2\epsilon}=\frac{-u_r-\sqrt{u_r^2+4\epsilon v_r}}{2\epsilon}
\label{3.16}
\end{equation}
and
\begin{equation}
\frac{-u_l+\sqrt{u_l^2+4\epsilon
v_l}}{2\epsilon}=\frac{-u_r+\sqrt{u_r^2+4\epsilon v_r}}{2\epsilon}
\label{3.17}
\end{equation}
respectively.

In order to identity the admissible solution,  the discontinuity
(\ref{3.14})  associating with $\lambda_1^\epsilon$ should satisfy
\begin{equation}
\sigma_1<\lambda_1^\epsilon(u_l, v_l)<\lambda_2^\epsilon(u_l, v_l),\
\ \ \ \lambda_1^\epsilon(u_r, v_r)<\sigma_1<\lambda_2^\epsilon(u_r,
v_r),\label{3.18}
\end{equation}
while the discontinuity (\ref{3.15}) associating with
$\lambda_2^\epsilon$ should satisfy
\begin{equation}
\lambda_1^\epsilon(u_l, v_l)<\sigma_2<\lambda_2^\epsilon(u_l, v_l),\
\ \ \ \lambda_1^\epsilon(u_r, v_r)<\lambda_2^\epsilon(u_r,
v_r)<\sigma_2.\label{3.19}
\end{equation}
Then one can check that both the inequality (\ref{3.18}) and
(\ref{3.19})  are equivalent to
\begin{equation}
u_r<u_l.\label{3.20}
\end{equation}

The discontinuity (\ref{3.14}) with (\ref{3.20}) is called as
1-shock and symbolized by $S_1$, and (\ref{3.15}) with (\ref{3.20})
is called as 2-shock and symbolized by $S_2$.

 For a given state $(u_l,v_l)$, all possible
states which can connect to $(u_l,v_l)$ on the right by a 1-shock
must be located on the straight line
\begin{equation}
{S_1(u_l,v_l)}:\ \ \ \ \
u-u_l=\Bigg(\displaystyle\frac{u_l-\sqrt{u_l^2+4\epsilon v_l}}{2
v_l}\ \Bigg)( v- v_l),\ \
 \ \ u<u_l,
\label{3.21}
\end{equation}
and all possible states which can  connect  to $(u_l,v_l)$ on the
right by a 2-shock must be located on the straight line
\begin{equation}
{S_2(u_l,v_l)}:\ \ \ \ \
 u- u_l=\Bigg(\displaystyle\frac{u_l+\sqrt{u_l^2+4\epsilon
   v_l}}{2 v_l}\Bigg)(v-v_l),\ \
 \ \ u<u_l.
\label{3.22}
\end{equation}

\vspace{0.2cm}

 Let us denote
$W_1(u_l,v_l)=R_1(u_l,v_l)\cup S_1(u_l,v_l)$ and
$W_2(u_l,v_l)=R_2(u_l,v_l)\cup S_2(u_l,v_l)$. Draw the line
$W_1(u_-,v_-)$ and $W_2(u_-,v_-)$ in the upper half  $(u,v)$-plane,
then the upper half  $(u,v)$-plane is divided into four regions (see
Fig.1). According to the right state $(u_+,v_+)$ in the different
regions, one can construct the unique global Riemann solution
connecting two constant states $(u_-,v_-)$ and $(u_+,v_+)$. To be
more exact, the Riemann solutions contain (i) a 1-rarefaction wave
and a 2-rarefaction wave when $(u_+,v_+)\in R_1R_2(u_-,v_-)$, (ii) a
1-rarefaction wave and a 2-shock
 when $(u_+,v_+)\in R_1S_2(u_-,v_-)$, (iii) a 1-shock and
a 2-rarefaction wave  when $(u_+,v_+)\in S_1R_2(u_-,v_-)$, (iv) a
1-shock and a 2-shock  when $(u_+,v_+)\in S_1S_2(u_-,v_-)$.

\begin{figure}[th]
\begin{center}
\includegraphics[width=10.0cm]{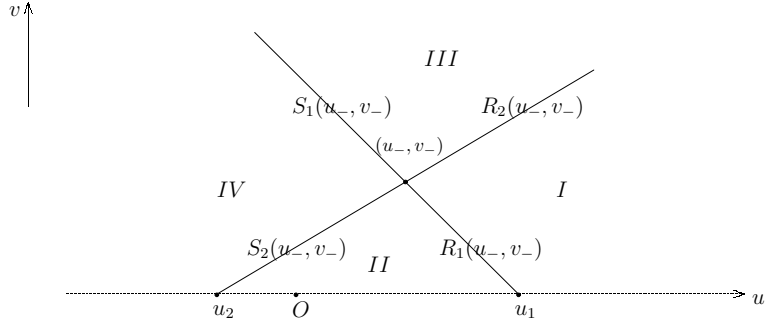}
\caption{Curves of elementary waves.  Here
$u_1=\frac{u_-+\sqrt{u_-^2+4\epsilon v_-}}{2}$ and
$u_2=\frac{u_--\sqrt{u_-^2+4\epsilon v_-}}{2}$}
\end{center}
\end{figure}

\vspace{0.2cm}

The conclusion can be stated in the following theorem.

\vspace{0.2cm}

\textsc{\textbf{Theorem 3.1.}} \textit{The Riemann problem for
$(\ref{DM})$ with initial data $(\ref{I-D})$  has  a unique
piecewise smooth solution consisting of waves of constant states,
shocks and rarefaction waves.}

\section{Limits of  Riemann solution to (\ref{DM}) for $u_ +<u_-,u_+/v_+<u_-/v_-$}

In this section, we study the limits of the Riemann solution
$\epsilon\rightarrow0^+$ when the initial data satisfy $u_ +<u_-,
u_+/v_+<u_-/v_-$. Especially, we pay more attention on the
phenomenon of concentration and the formation of delta-shocks  in
the limit.

For  $u_ +<u_-,u_+/v_+<u_-/v_-$, there must exist $\epsilon_0>0$
such that the Riemann solution just consists of two shocks   for any
$\epsilon<\epsilon_0$. In fact, since all states $(u,v)$ connected
with $(u_-,v_-)$ by $S_1$ and $S_2$ satisfy
$$
u-u_-=\Bigg(\displaystyle\frac{u_--\sqrt{u_-^2+4\epsilon v_-}}{2
v_-}\ \Bigg)( v- v_-),\ \
 \ \ u<u_-,\ \
 \ \ v>v_-,
$$
and
$$
 u- u_-=\Bigg(\displaystyle\frac{u_-+\sqrt{u_-^2+4\epsilon
   v_-}}{2 v_-}\Bigg)(v-v_-),\ \
 \ \ u<u_-,\ \
 \ \ v<v_-,
$$
respectively, then if $v_+=v_-$, $\epsilon_0$ may be taken any real
positive number, otherwise, we have  the conclusion by taking
$$
\epsilon_0=\frac{\Big(2v_-\Big(\frac{u_+-u_-}{v_+-v_-}\Big)-u_-\Big)^2-u_-^2}{4v_-}
        =\frac{(u_+-u_-)(v_-u_+-v_+u_-)}{(v_+-v_-)^2}.
$$

\vspace{0.2cm}

For fixed $\epsilon<\epsilon_0$, let $ {U}^\epsilon(\xi)$ denote the
two-shock Riemann solution for  (\ref{DM}) and (\ref{I-D})
constructed in Section 3
\begin{equation}
 {U}^\epsilon(\xi)=(u^{\epsilon},v^{\epsilon})(\xi)=\left\{\begin{array}{ll}
 (u_-,v_-),&\xi<\sigma_1^{\epsilon}, \\[1mm]
 (u_*^{\epsilon},v_*^{\epsilon}), &\sigma_1^{\epsilon}<\xi<\sigma_2^{\epsilon},\\[1mm]
 (u_+,v_+), &\xi>\sigma_2^{\epsilon},
 \end{array}\right.\label{4.1}
\end{equation}
where $(u_-,v_-)$ and $(u_*^{\epsilon},v_*^{\epsilon})$ are
connected by a shock $S_1$ with speed $\sigma_1^{\epsilon}$, and
$(u_*^{\epsilon},v_*^{\epsilon})$ and $(u_+,v_+)$ are connected by a
shock $S_2$ with speed $\sigma_2^{\epsilon}$:
\begin{equation}
S_1: \left\{\begin{array}{l}
 \sigma_1^\epsilon=u_*^{\epsilon}+\displaystyle\frac{u_--\sqrt{u_-^2+4\epsilon v_-}}{2},\cr\noalign {\vskip2truemm}
 u_*^{\epsilon}-u_-=\Bigg(\displaystyle\frac{u_--\sqrt{u_-^2+4\epsilon v_-}}{2 v_-}\Bigg)(v_*^{\epsilon}-
 v_-),
\end{array}\ \ \ \ \ v_*^{\epsilon}>v_-,
\right.\label{4.2}
\end{equation}
\begin{equation}
S_2: \left\{\begin{array}{l}
 \sigma_2^\epsilon=u_++\displaystyle\frac{u_*^{\epsilon}+\sqrt{(u_*^{\epsilon})^2+4\epsilon v_*^{\epsilon}}}{2},\cr\noalign{\vskip4truemm}
  u_+- u_*^{\epsilon}=\Bigg(\displaystyle\frac{u_*^{\epsilon}+\sqrt{(u_*^{\epsilon})^2+4\epsilon v_*^{\epsilon}}}{2 v_*^{\epsilon}}\Bigg)(v_+-v_*^{\epsilon}),
\end{array}\ \ \ \ \ v_*^{\epsilon}>v_+.
\right.\label{4.3}
\end{equation}
Here
\begin{equation}
\frac{u_--\sqrt{u_-^2+4\epsilon
v_-}}{2v_-}=\frac{u_*^{\epsilon}-\sqrt{(u_*^{\epsilon})^2+4\epsilon
v_*^{\epsilon}}}{2v_*^{\epsilon}}\label{4.4}
\end{equation}
and
\begin{equation}
\frac{u_*^{\epsilon}+\sqrt{(u_*^{\epsilon})^2+4\epsilon
v_*^{\epsilon}}}{2v_*^{\epsilon}}=\frac{u_++\sqrt{u_+^2+4\epsilon
v_+}}{2v_+}.\label{4.4}
\end{equation}

\vspace{0.2cm}

The following Lemmas 4.1-4.2 show  the limit behaviors of the
 states between two shocks.

\vspace{0.2cm}

\textsc{\textbf{Lemma 4.1.}}
$$
\lim\limits_{\epsilon\rightarrow0^+}v_*^\epsilon=
\left\{\begin{array}{ll}
(u_-/u_+)v_+,& for \ \ \ u_->u_+>0,\\
(u_+/u_-)v_-,& for \ \ \ 0>u_->u_+,\\
+\infty,& for \ \ \ u_-\geq0\geq u_+.
\end{array}\right.
$$

\vspace{0.2cm}

\emph{Proof.} Based on (\ref{4.2}) and (\ref{4.3}), $v_*^\epsilon$
can be expressed as
$$
 u_-+\Bigg(\displaystyle\frac{u_--\sqrt{u_-^2+4\epsilon v_-}}{2 v_-}\Bigg)(v_*^{\epsilon}-v_-)=
u_+-\Bigg(\frac{u_++\sqrt{u_+^2+4\epsilon
v_+}}{2v_+}\Bigg)(v_+-v_*^{\epsilon}).
$$
Solving this equation gives
\begin{equation}
v_*^{\epsilon}=\frac{u_+-u_--\Big(\frac{u_++\sqrt{u_+^2+4\epsilon
v_+}}{2v_+}\Big)v_++\Big(\frac{u_--\sqrt{u_-^2+4\epsilon v_-}}{2
v_-}\Big)v_-} {\frac{u_--\sqrt{u_-^2+4\epsilon v_-}}{2
v_-}-\frac{u_++\sqrt{u_+^2+4\epsilon v_+}}{2v_+}}.\label{4.6}
\end{equation}
Taking the limit $\epsilon\rightarrow0^+$ will lead to the
conclusions. The proof is finished.  \hspace{1cm}$\blacksquare$

\vspace{0.2cm}

\textsc{\textbf{Lemma 4.2.}}
$$
\lim\limits_{\epsilon\rightarrow0^+}u_*^\epsilon=
\left\{\begin{array}{ll}
u_-,& for \ \ \ u_->u_+>0,\\
u_+,& for \ \ \ 0>u_->u_+,\\
u_-+u_+,& for \ \ \ u_-\geq0\geq u_+.
\end{array}\right.
$$

\vspace{0.2cm}

\emph{Proof.} From the second equation in (\ref{4.2}), we have
$$
 u_*^{\epsilon}=u_-+\Bigg(\displaystyle\frac{u_--\sqrt{u_-^2+4\epsilon v_-}}{2 v_-}\Bigg)(v_*^{\epsilon}-
 v_-).
$$
For the cases $u_->u_+>0$ and $0>u_->u_+$, the conclusions are
obvious because of the Lemma 4.2. For the case $u_-\geq0\geq u_+$,
due to
$$
\frac{\frac{u_++\sqrt{u_+^2+4\epsilon
v_+}}{2v_+}}{\frac{u_--\sqrt{u_-^2+4\epsilon v_-}}{2 v_-}}=
\displaystyle\frac{u_++\sqrt{u_+^2+4\epsilon
 v_+}}{u_--\sqrt{u_-^2+4\epsilon v_-}}\cdot\frac{v_-}{v_+}=\displaystyle\frac{u_-+\sqrt{u_-^2+4\epsilon v_-}}{u_+-\sqrt{u_+^2+4\epsilon v_+}}\rightarrow\frac{u_-}{u_+}\ \ \ \ \ as\ \epsilon\rightarrow0^+,
$$
we have
$$
\begin{array}{ll}
\Bigg(\frac{u_--\sqrt{u_-^2+4\epsilon v_-}}{2
v_-}\Bigg)v_*^{\epsilon}&= \Bigg(\frac{u_--\sqrt{u_-^2+4\epsilon
v_-}}{2
v_-}\Bigg)\cdot\frac{u_+-u_--\Big(\frac{u_++\sqrt{u_+^2+4\epsilon
v_+}}{2v_+}\Big)v_++\Big(\frac{u_--\sqrt{u_-^2+4\epsilon v_-}}{2
v_-}\Big)v_-}
{\frac{u_--\sqrt{u_-^2+4\epsilon v_-}}{2 v_-}-\frac{u_++\sqrt{u_+^2+4\epsilon v_+}}{2v_+}}\\[6mm]
 &=\frac{u_+-u_--\Big(\frac{u_++\sqrt{u_+^2+4\epsilon v_+}}{2v_+}\Big)v_++\Big(\frac{u_--\sqrt{u_-^2+4\epsilon v_-}}{2 v_-}\Big)v_-}
{1-\frac{\frac{u_++\sqrt{u_+^2+4\epsilon
v_+}}{2v_+}}{\frac{u_--\sqrt{u_-^2+4\epsilon v_-}}{2
v_-}}}\rightarrow u_+ \ \ \ \ \ as\ \epsilon\rightarrow0^+,
\end{array}
$$
which gives the conclusion. The proof is finished.
\hspace{1cm}$\blacksquare$

\vspace{0.2cm}

The following Lemma 4.3 shows the limit behaviors of the speeds of
two shocks.

\vspace{0.2cm}

\textsc{\textbf{Lemma 4.3.}}
$$
\lim\limits_{\epsilon\rightarrow0^+}(\sigma_1^\epsilon,\sigma_2^\epsilon)=
\left\{\begin{array}{ll}
(u_-,u_-+u_+),& for \ \ \ u_->u_+>0,\\
(u_-+u_+,u_+),& for \ \ \ 0>u_->u_+,\\
(u_-+u_+,u_-+u_+),& for \ \ \ u_-\geq0\geq u_+.
\end{array}\right.
$$

\vspace{0.2cm}

\emph{Proof.} For $u_->u_+>0$,
$$
\begin{array}{l}
\lim\limits_{\epsilon\rightarrow0^+}\sigma_1^\epsilon=\lim\limits_{\epsilon\rightarrow0^+}\Bigg(u_*^{\epsilon}+\frac{u_--\sqrt{u_-^2+4\epsilon
 v_-}}{2}\Bigg)
  =\lim\limits_{\epsilon\rightarrow0^+}u_*^{\epsilon}
  =u_-,
 \end{array}
$$
and
$$
\begin{array}{ll}
 \lim\limits_{\epsilon\rightarrow0^+}\sigma_2^\epsilon&=\lim\limits_{\epsilon\rightarrow0^+}\Bigg(u_++\displaystyle\frac{u_*^{\epsilon}+\sqrt{(u_*^{\epsilon})^2+4\epsilon v_*^{\epsilon}}}{2}\Bigg)\cr\noalign{\vskip2truemm}
 &=\lim\limits_{\epsilon\rightarrow0^+}\Bigg(u_++\Bigg(\displaystyle\frac{u_++\sqrt{u_+^2+4\epsilon v_+}}{2v_+}\Bigg)v_*^\epsilon\Bigg)\cr\noalign{\vskip4truemm}
 &=u_-+u_+.
\end{array}
$$
For $0>u_->u_+$, the conclusions can be proved in a similar way. For
$u_-\geq0\geq u_+$,
$$
\lim\limits_{\epsilon\rightarrow0^+}\sigma_1^\epsilon=\lim\limits_{\epsilon\rightarrow0^+}u_*^{\epsilon}
=u_-+u_+
$$
and
$$
\begin{array}{ll}
 \lim\limits_{\epsilon\rightarrow0^+}\sigma_2^\epsilon&=\lim\limits_{\epsilon\rightarrow0^+}\Bigg(u_++\displaystyle\frac{u_*^{\epsilon}+\sqrt{(u_*^{\epsilon})^2+4\epsilon v_*^{\epsilon}}}{2}\Bigg)\cr\noalign{\vskip2truemm}
 &=\lim\limits_{\epsilon\rightarrow0^+}\Bigg(u_++\Bigg(\displaystyle\frac{u_++\sqrt{u_+^2+4\epsilon v_+}}{2v_+}\Bigg)v_*^\epsilon\Bigg)\cr\noalign{\vskip4truemm}
 &=\lim\limits_{\epsilon\rightarrow0^+}\Bigg(u_++\Bigg(\displaystyle\frac{u_++\sqrt{u_+^2+4\epsilon v_+}}{2v_+}\Bigg)v_+-(u_+-u_*^\epsilon)\Bigg)\cr\noalign{\vskip4truemm}
 &=u_-+u_+,
\end{array}
$$
where we have used
 \begin{equation}
\begin{array}{l}
 \Bigg( \frac{u_++\sqrt{u_+^2+4\epsilon
v_+}}{2v_+}\Bigg)v_*^\epsilon=\Bigg( \frac{u_++\sqrt{u_+^2+4\epsilon
v_+}}{2v_+}\Bigg)v_+-(u_+-u_*^\epsilon)
\end{array}\label{4.7}
\end{equation}
obtaining from the second equality of (\ref{4.3}). The proof is
finished. \hspace{1cm}$\blacksquare$

\vspace{0.2cm}

Let $U^0(\xi)=\lim\limits_{\epsilon\rightarrow0^+}U^\epsilon(\xi)$.
Then when $u_->u_+>0, u_+/v_+<u_-/v_-$,
$$
U^0(\xi)=\left\{\begin{array}{ll}
(u_-,v_-),&\xi<\sigma_1,\\
(u_-,v_*),&\sigma_1<\xi<\sigma_2,\\
(u_+,v_+),&\xi>\sigma_2,
\end{array}\right.
$$
where $\sigma_1=u_-, \sigma_2=u_-+u_+$ and  $v_*=(u_-/u_+)v_-$. When
$0>u_->u_+, u_+/v_+<u_-/v_-$,
$$
U^0(\xi)=\left\{\begin{array}{ll}
(u_-,v_-),&\xi<\sigma_1,\\
(u_+,v_*),&\sigma_1<\xi<\sigma_2,\\
(u_+,v_+),&\xi>\sigma_2,
\end{array}\right.
$$
$\sigma_1=u_-+u_+, \sigma_2=u_+$ and  $v_*=(u_+/u_-)v_+$. It can be
seen that $U^0(\xi)$ coincides with the Riemann solution for
(\ref{Bur}) constructed in Section 2.

\vspace{0.2cm}

For the case $u_-\geq0\geq  u_+$, it has been shown that two shocks
will coincide at $\xi=u_-+u_+:=\sigma$ as $\epsilon\rightarrow0^+$.
Furthermore, for  the component $u^\epsilon(\xi)$, it has  been
shown that
\begin{equation}
\lim\limits_{\epsilon\rightarrow0^+}u^\epsilon(\xi)=\left\{\begin{array}{ll}
u_-,&\xi<\sigma,\\
u_-+u_+,&\xi=\sigma,\\
u_+,&\xi>\sigma.
\end{array}\right.\label{4.8}
\end{equation}
For  the component $v^\epsilon(\xi)$, we have proven that the
intermediate state $v^\epsilon_*$ becomes infinity  as
$\epsilon\rightarrow0^+$. Further, we have

\vspace{0.2cm}

\textsc{\textbf{Lemma 4.4.}}
$$
\lim\limits_{\epsilon\rightarrow0^+}
(\sigma_2^\epsilon-\sigma_1^\epsilon)v_*^\epsilon=u_-v_+-u_+v_-.
$$

\vspace{0.2cm}

\emph{Proof.} With
$$
\displaystyle\frac{u_*^{\epsilon}+\sqrt{(u_*^{\epsilon})^2+4\epsilon
v_*^{\epsilon}}}{2}=
\displaystyle\frac{u_*^{\epsilon}+\sqrt{(u_*^{\epsilon})^2+4\epsilon
v_*^{\epsilon}}}{2v_*^{\epsilon}}v_*^{\epsilon}
=\displaystyle\frac{u_++\sqrt{u_+^2+4\epsilon
v_+}}{2v_+}v_*^{\epsilon}
$$
and (\ref{4.7}), it follows
$$
\begin{array}{ll}
\sigma_2^\epsilon-\sigma_1^\epsilon
&=u_++\displaystyle\frac{u_*^{\epsilon}+\sqrt{(u_*^{\epsilon})^2+4\epsilon
v_*^{\epsilon}}}{2}-
u_*^{\epsilon}-\displaystyle\frac{u_--\sqrt{u_-^2+4\epsilon v_-}}{2}\\[4mm]
&=\displaystyle\frac{u_++\sqrt{u_+^2+4\epsilon
v_+}}{2}-\displaystyle\frac{u_--\sqrt{u_-^2+4\epsilon v_-}}{2}.
\end{array}
$$
Then
$$
\begin{array}{ll}
\lim\limits_{\epsilon\rightarrow0^+}(\sigma_2^\epsilon-\sigma_1^\epsilon)v_*^\epsilon
 &=\lim\limits_{\epsilon\rightarrow0^+}\Big(\frac{u_++\sqrt{u_+^2+4\epsilon v_+}}{2}-\frac{u_--\sqrt{u_-^2+4\epsilon v_-}}{2}\Big)
               \cdot\frac{\frac{u_+-\sqrt{u_+^2+4\epsilon v_+}}{2}-\frac{u_-+\sqrt{u_-^2+4\epsilon v_-}}{2}}{\frac{u_--\sqrt{u_-^2+4\epsilon v_-}}{2v_-}-\frac{u_++\sqrt{u_+^2+4\epsilon v_+}}{2v_+}}\\[6mm]
  &=\lim\limits_{\epsilon\rightarrow0^+}\Big(\frac{u_+-\sqrt{u_+^2+4\epsilon v_+}}{2}-\frac{u_-+\sqrt{u_-^2+4\epsilon v_-}}{2}\Big)\cdot
 \frac{(u_++\sqrt{u_+^2+4\epsilon v_+})-(u_--\sqrt{u_-^2+4\epsilon v_-})}{\frac{u_--\sqrt{u_-^2+4\epsilon v_-}}{v_-}-\frac{u_++\sqrt{u_+^2+4\epsilon v_+}}{v_+}}\\[6mm]
 &=\lim\limits_{\epsilon\rightarrow0^+}\Big(\frac{u_+-\sqrt{u_+^2+4\epsilon v_+}}{2}-\frac{u_-+\sqrt{u_-^2+4\epsilon v_-}}{2}\Big)\cdot
 \frac{\frac{u_++\sqrt{u_+^2+4\epsilon v_+}}{u_--\sqrt{u_-^2+4\epsilon v_-}}-1}{\frac{1}{v_-}-\frac{1}{v_+}\cdot\frac{u_++\sqrt{u_+^2+4\epsilon v_+}}{u_--\sqrt{u_-^2+4\epsilon v_-}}}\\[4mm]
  &=u_-v_+-u_+v_- .
\end{array}
$$
The proof is finished.  \hspace{1cm}$\blacksquare$

\vspace{0.4cm}

Take $\phi(\xi)\in C^\infty_0(-\infty,+\infty)$ such that
$\phi(\xi)\equiv\phi(\sigma)$ for $\xi$ in a neighborhood $\Omega$
of $\xi=\sigma$ ($\phi$ is called a sloping test function
\cite{D-M,G-T}). Assume when $\epsilon<\epsilon_0$, it holds
$\sigma_1^\epsilon\in \Omega$ and $\sigma_2^\epsilon\in \Omega$. It
is well known that the solution (\ref{4.1}) satisfies weak
formulation
 \begin{equation}
-\displaystyle\int^{+\infty}_{-\infty}v^{\epsilon}(u^{\epsilon}-\xi)\phi'd\xi+\displaystyle\int^{+\infty}_{-\infty}v^{\epsilon}\phi
d\xi=0.\label{4.9}
\end{equation}
Since
$$
\displaystyle\int^{+\infty}_{-\infty}v^{\epsilon}(u^{\epsilon}-\xi)\phi'd\xi=\bigg(\displaystyle\int^{\sigma_1^{\epsilon}}_{-\infty}
+\displaystyle\int^{+\infty}_{\sigma_2^{\epsilon}}\bigg)v^{\epsilon}(u^{\epsilon}-\xi)\phi'd\xi,
$$
we have
$$
\begin{array}{ll}
 \lim\limits_{\epsilon\rightarrow0^+}\displaystyle\int^{+\infty}_{-\infty}v^{\epsilon}(u^{\epsilon}-\xi)\phi'd\xi
 &=\lim\limits_{\epsilon\rightarrow0^+}\displaystyle\int^{\sigma_1^{\epsilon}}_{-\infty}v_-(u_--\xi)\phi'd\xi+\lim\limits_{\epsilon\rightarrow0^+}\displaystyle\int_{\sigma_2^{\epsilon}}^{+\infty}v_+(u_+-\xi)\phi'd\xi\\[6mm]
 &=\Big(u_-v_+-u_+v_-\Big)\phi(\sigma)+\displaystyle\int^{+\infty}_{-\infty}H(\xi-\sigma)\phi  d\xi,
\end{array}
$$
where
$$
H(x)=\left\{\begin{array}{ll}
 v_-,&x<0,\\[1mm]
 v_+,&x>0.
 \end{array}\right.
$$
Returning to (\ref{4.9}),we get
 \begin{equation}
\lim\limits_{\epsilon\rightarrow0^+}\displaystyle\int^{+\infty}_{-\infty}\Big(v^{\epsilon}-H(\xi-\sigma)\Big)\phi
d\xi =(u_-v_+-u_+v_-)\phi(\sigma)\label{4.10}
\end{equation}
for all sloping test functions $\phi\in
C^\infty_0(-\infty,+\infty)$.

For an arbitrary test function $\varphi(\xi)\in
C^\infty_0(-\infty,+\infty)$, we take a sloping test function $\phi$
such that $\phi(\sigma)=\varphi(\sigma)$ and
$$
\max\limits_{\xi\in (-\infty,+\infty)}|\phi-\varphi|<\mu.
$$
We have
$$
\begin{array}{l}
\lim\limits_{\epsilon\rightarrow0^+}\displaystyle\int^{+\infty}_{-\infty}\Big(v^{\epsilon}-H(\xi-\sigma)\Big)\varphi d\xi\\[3mm]
 =\lim\limits_{\epsilon\rightarrow0^+}\displaystyle\int^{+\infty}_{-\infty}\Big(v^{\epsilon}-H(\xi-\sigma)\Big)\phi d\xi+\lim\limits_{\epsilon\rightarrow0^+}\displaystyle\int^{+\infty}_{-\infty}\Big(v^{\epsilon}-H(\xi-\sigma)\Big)(\varphi-\phi) d\xi.
\end{array}
$$
The first limit on the right side
$$
\begin{array}{ll}
\lim\limits_{\epsilon\rightarrow0^+}\displaystyle\int^{+\infty}_{-\infty}\Big(v^{\epsilon}-H(\xi-\sigma)\Big)\phi
d\xi
 &=(u_-v_+-u_+v_-)\phi(\sigma)\\[3mm]
&=(u_-v_+-u_+v_-)\varphi(\sigma).
\end{array}
$$
The second limit on the right side
$$
\begin{array}{ll}
\displaystyle\int^{+\infty}_{-\infty}\Big(v^{\epsilon}_*-H(\xi-\sigma)\Big)(\varphi-\phi)
d\xi
&=\displaystyle\int^{\sigma_2^\epsilon}_{\sigma_1^\epsilon}\Big(v^{\epsilon}_*-H(\xi-\sigma)\Big)(\varphi-\phi) d\xi\\[3mm]
&=\displaystyle\int^{\sigma_2^\epsilon}_{\sigma_1^\epsilon}v^{\epsilon}(\varphi-\phi)
d\xi+
\displaystyle\int^{\sigma_2^\epsilon}_{\sigma_1^\epsilon}H(\xi-\sigma)(\varphi-\phi) d\xi,\\[3mm]
\end{array}
$$
which converges to 0 by sending $\mu\rightarrow0$ and recalling
Lemma 4.5. Thus we have that
 \begin{equation}
\lim\limits_{\epsilon\rightarrow0^+}\displaystyle\int^{+\infty}_{-\infty}\Big(v^{\epsilon}-H(\xi-\sigma)\Big)\varphi
d\xi =(u_-v_+-u_+v_-)\varphi(\sigma)\label{4.11}
\end{equation}
for all test functions $\varphi\in C^\infty_0(-\infty,+\infty)$.

 Let $\psi(x,t)\in C^\infty_0((-\infty,+\infty)\times[0,+\infty))$ be
a smooth test function, and let $\widetilde{\psi}(\xi,t):=\psi(\xi
t,t)$. Then it follows
$$
\begin{array}{ll}
\lim\limits_{\epsilon\rightarrow0^+}\displaystyle\int^{+\infty}_{0}\displaystyle\int^{+\infty}_{-\infty}v^{\epsilon}(x/t)\psi(x,t)dxdt
 &=\lim\limits_{\epsilon\rightarrow0^+}\displaystyle\int^{+\infty}_{0}\displaystyle\int^{+\infty}_{-\infty}v^{\epsilon}(\xi)\psi(\xi t,t)d(\xi t)dt\\[5mm]
 &=\lim\limits_{\epsilon\rightarrow0^+}\displaystyle\int^{+\infty}_{0}t\bigg(\displaystyle\int^{+\infty}_{-\infty}v^{\epsilon}(\xi)\widetilde{\psi}(\xi,t)d\xi\bigg)dt
\end{array}
$$
and from (\ref{4.11})
$$
\begin{array}{l}
 \lim\limits_{\epsilon\rightarrow0^+}\displaystyle\int^{+\infty}_{-\infty}v^{\epsilon}(\xi)\widetilde{\psi}(\xi,t)d\xi\\[5mm]
 =\displaystyle\int^{+\infty}_{-\infty}H (\xi-\sigma)\widetilde{\psi}(\xi,t) d\xi
     +\big(u_-v_+-u_+v_-\big)\widetilde{\psi}(\sigma,t)\\[5mm]
 =t^{-1}\displaystyle\int^{+\infty}_{-\infty}H (x-\sigma t)\psi(x,t)dx+\big(u_-v_+-u_+v_-\big)\psi(\sigma t,t).
\end{array}
$$
Combining the two relations above yields
$$
\begin{array}{l}
\lim\limits_{\epsilon\rightarrow0^+}\displaystyle\int^{+\infty}_{0}\displaystyle\int^{+\infty}_{-\infty}v^{\epsilon}(x/t)\psi(x,t)dxdt\\[5mm]
 =\displaystyle\int^{+\infty}_{0}\displaystyle\int^{+\infty}_{-\infty}
     H (x-\sigma t)\psi(x,t)dxdt+\displaystyle\int^{+\infty}_{0}\big(u_-v_+-u_+v_-\big)t\psi(\sigma  t,t)dt.
\end{array}
$$
The last term, by the definition,
$$
 \displaystyle\int^{+\infty}_{0}\big(u_-v_+-u_+v_-\big)t\psi(\sigma t,t)dt=\Big<w(t)\delta_{x=\sigma t},\psi(x,t)\Big>
$$
with
$$
w(t)=\big(u_-v_+-u_+v_-\big)t.
$$

\vspace{0.2cm}

Thus we obtain the following Theorem.

\vspace{0.2cm}

\textsc{\textbf{Theorem 4.5.}} \emph{Let $u_-\geq0\geq u_+$, and
$(u^\epsilon(x,t),v^\epsilon(x,t))$ is the two-shock solution to
$(\ref{DM})$ and $(\ref{I-D})$. Then
$(u^\epsilon(x,t),v^\epsilon(x,t))$ converges in the sense of
distributions. Denote the limit functions $U^0(x,t)$, then
$$
U^0(x,t)=\left\{\begin{array}{ll}
(u_-,v_-),&x<\sigma t,\\
\big(u_-+u_+,w(t)\delta(x-\sigma t)\big),&x=\sigma t,\\
(u_+,v_+),&x>\sigma t,
\end{array}\right.
$$
where $w(t)=(u_-v_+-u_+v_-)t$  and $\sigma=u_-+u_+$, which  is just
the delta-shock Riemann solution of (\ref{Bur}) with the same
initial data.}

\section{Limits  of Riemann  solution  to (\ref{DM})  for $u_ +>u_-, u_+/v_+>u_-/v_-$}

In this section, we study the limits  of the Riemann solution as
$\epsilon\rightarrow0^+$ when the initial data satisfy $u_ +>u_-,
u_+/v_+>u_-/v_-$. At this moment,  there must exist $\epsilon_0>0$
such that the Riemann solution just consists of two rarefaction
waves  for any $\epsilon<\epsilon_0$.

For fixed $\epsilon<\epsilon_0$, let $ {U}^\epsilon(\xi)$ denote the
two-rarefaction-wave Riemann solution for (\ref{DM})  and
(\ref{I-D}) constructed in Section 3
\begin{equation}
 {U}^\epsilon(\xi)=(u^{\epsilon},v^{\epsilon})(\xi)=\left\{\begin{array}{cl}
 (u_-,v_-),&-\infty<\xi\leq \lambda_1(u_-,v_-), \\[1mm]
 R_1,&\lambda_1(u_-,v_-)\leq\xi\leq\lambda_1(u_*^\epsilon,v_*^\epsilon), \\[1mm]
 (u_*^{\epsilon},v_*^{\epsilon}), &\lambda_1(u_*^\epsilon,v_*^\epsilon)\leq\xi\leq\lambda_2(u_*^\epsilon,v_*^\epsilon),\\[1mm]
 R_2,&\lambda_2(u_*^\epsilon,v_*^\epsilon)\leq\xi\leq\lambda_2(u_+,v_+), \\[1mm]
 (u_+,v_+), &\lambda_2(u_+,v_+)\leq\xi<+\infty,
 \end{array}\right.\label{5.1}
\end{equation}
where
\begin{equation}
R_1: \left\{\begin{array}{l}
 \xi=\lambda_1^\epsilon(u,v)=u+\displaystyle\frac{u-\sqrt{u^2+4\epsilon v}}{2},\cr\noalign {\vskip4truemm}
  u-u_-=\Bigg(\displaystyle\frac{u_--\sqrt{u_-^2+4\epsilon v_-}}{2 v_-}\ \Bigg)( v- v_-),
\end{array}v_-\geq v\geq v_*^{\epsilon},\ u_*^{\epsilon}\geq u\geq u_-,
\right.\label{5.2}
\end{equation}
and
\begin{equation}
R_2: \left\{\begin{array}{l}
 \xi=\lambda_2^\epsilon(u,v)=u+\displaystyle\frac{u+\sqrt{u^2+4\epsilon v}}{2},\cr\noalign {\vskip2truemm}
   u- u_+=\Bigg(\displaystyle\frac{u_++\sqrt{u_+^2+4\epsilon
   v_+}}{2 v_+}\Bigg)(v-v_+),
\end{array}v_+\geq v\geq v_*^{\epsilon},\ u_+\geq u\geq u_*^{\epsilon}.
\right.\label{5.3}
\end{equation}

\vspace{0.2cm}

The follow lemmas describe  the limit behaviors of the intermediate
state $(u_*^{\epsilon},v_*^{\epsilon})$ between two rarefaction
waves.

\vspace{0.2cm}

\textsc{\textbf{Lemma 5.1.}}
$$
\lim\limits_{\epsilon\rightarrow0^+}v_*^\epsilon=
\left\{\begin{array}{ll}
(u_-/u_+)v_+,& for \ \ \ u_+>u_->0,\\
(u_+/u_-)v_-,& for \ \ \ 0>u_+>u_-,\\
0,& for \ \ \ u_+>0>u_-.
\end{array}\right.
$$

\vspace{0.2cm}

\emph{Proof.} From (\ref{5.2}) and (\ref{5.3}), it follows
$$
  v_*^\epsilon=\frac{u_+-u_--\frac{u_++\sqrt{u_+^2+4\epsilon v_+}}{2
v_+} v_++\frac{u_--\sqrt{u_-^2+4\epsilon v_-}}{2 v_-}
v_-}{\frac{u_--\sqrt{u_-^2+4\epsilon v_-}}{2 v_-}
-\frac{u_++\sqrt{u_+^2+4\epsilon v_+}}{2 v_+}}.
$$
Then the conclusions can be obtained directly by taking the limit
$\epsilon\rightarrow0^+$. The proof is finished.
\hspace{1cm}$\blacksquare$

\vspace{0.2cm}

\textsc{\textbf{Lemma 5.2.}}
$$
\lim\limits_{\epsilon\rightarrow0^+}u_*^\epsilon=
\left\{\begin{array}{ll}
u_-,& for \ \ \ u_+>u_->0,\\
u_+,& for \ \ \ 0>u_+>u_-,\\
0,& for \ \ \ u_+\geq0\geq u_-.
\end{array}\right.
$$

\vspace{0.2cm}

\emph{Proof.} From (\ref{5.2}), we have
$$
  u=u_-+\Bigg(\displaystyle\frac{u_--\sqrt{u_-^2+4\epsilon v_-}}{2 v_-}\ \Bigg)(v_*^\epsilon- v_-).
$$
With the Lemma 5.1, we easily get  the conclusions. The proof is
complete.  \hspace{1cm}$\blacksquare$

\vspace{0.2cm}

Besides, as $\epsilon\rightarrow0^+$, when $u_+>u_->0$, the
rarefaction wave $R_1$ tends to
\begin{equation}
  \begin{array}{l}
 \xi=u=u_-,
\end{array} \label{5.4}
\end{equation}
and the rarefaction wave $R_2$ tends to
\begin{equation}
 \begin{array}{l}
 \xi=2u,\ \
   u/v=u_+/v_+.
\end{array} \label{5.5}
\end{equation}
When $0>u_+>u_-$, the rarefaction wave $R_1$ tends to
\begin{equation}
 \begin{array}{l}
 \xi=2u,\ \ \
  u/v=u_-/v_-,
\end{array}
 \label{5.6}
\end{equation}
and  the rarefaction wave $R_2$ tends to
\begin{equation}
 \begin{array}{l}
 \xi=u=u_+,
\end{array}
 \label{5.7}
\end{equation}
When $u_+\geq0\geq u_-$, the rarefaction wave $R_1$ tends to
(\ref{5.6}), and the rarefaction wave $R_2$ tends to (\ref{5.5}).

 \vspace{0.2cm}

In conclusion,   when $u_ +>u_-, u_+/v_+>u_-/v_-$, the limits of
Riemann solution of (\ref{DM}) are just the solutions of (\ref{Bur})
with the same initial data.

\vspace{0.3cm}

In the above two sections, we have proven that when $u_ +<u_-,
u_+/v_+<u_-/v_-$ and $u_ +>u_-, u_+/v_+>u_-/v_-$, the solutions to
the Riemann problem for (\ref{DM}) just are the solutions to the
Riemann problem for (\ref{Bur}) with the same initial data. The same
conclusions are true for the rest two cases $u_ +>u_-,
u_+/v_+<u_-/v_-$ and $u_ +<u_-, u_+/v_+>u_-/v_-$, and we omit the
discussions.

\section {Process of concentration: Numerical simulations}

To understand the phenomenon of concentration and the process of
formation of delta-shocks in the Riemann solutions to (\ref{DM}) as
the flux $\epsilon v$ vanishes, in this section, we present some
representative numerical results, obtained by  employing the
Nessyahu-Tadmor scheme \cite{N-T,Jiang} with $500$ cells and CFL =
0.475.  We take the  initial data as follows
$$
(u,v)(x,t=0)=\left\{\begin{array}{ll}
  (1,1),&x<0,\\
  (-1,1.5), &x>0.
 \end{array}\right.
$$
The numerical simulations for  different choices of $\epsilon$ are
presented in Figs.2-5.

\begin{figure}[th]
\begin{center}
\includegraphics[width=2.5in]{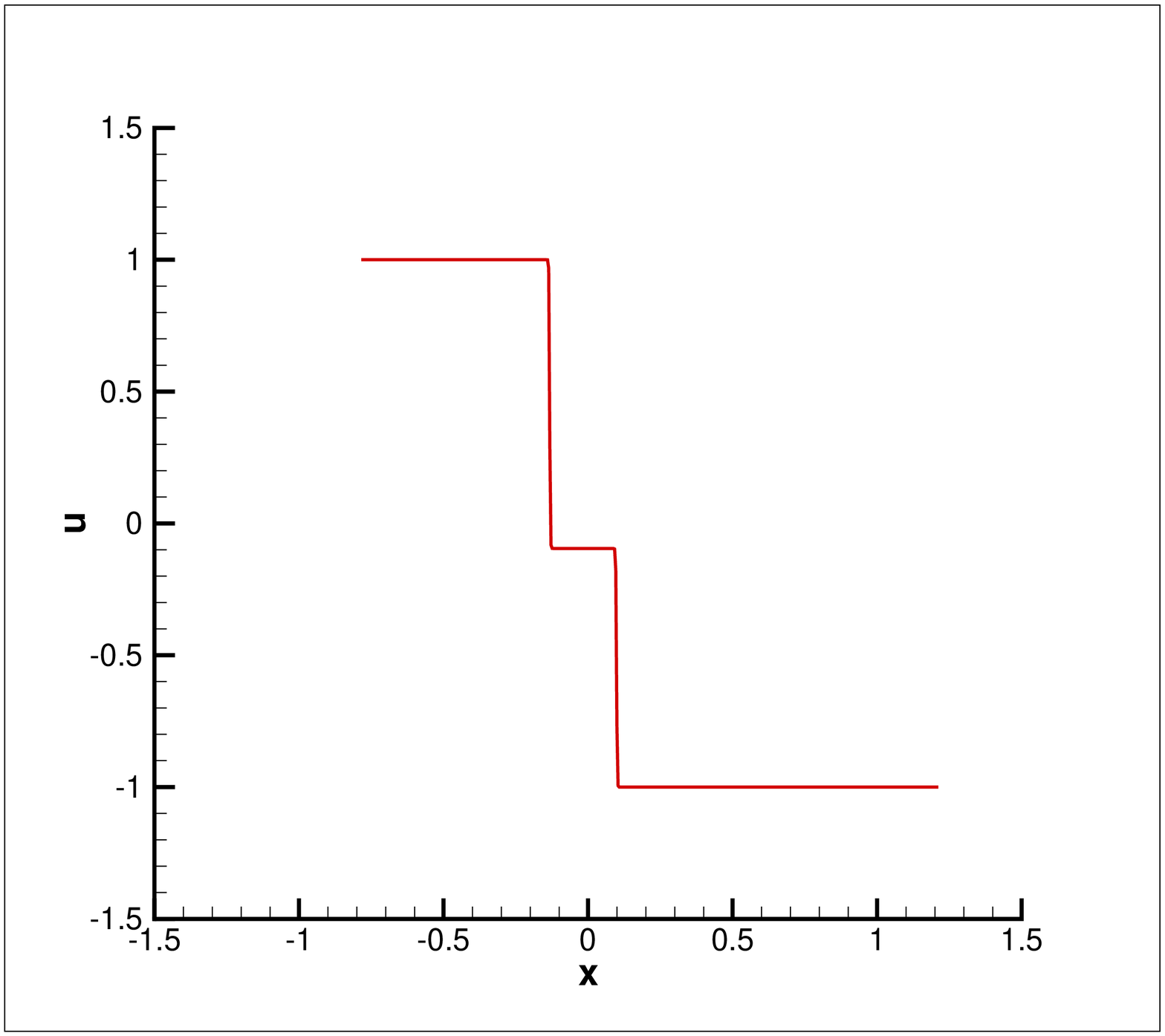}\hspace{0.50cm}\includegraphics[width=2.5in]{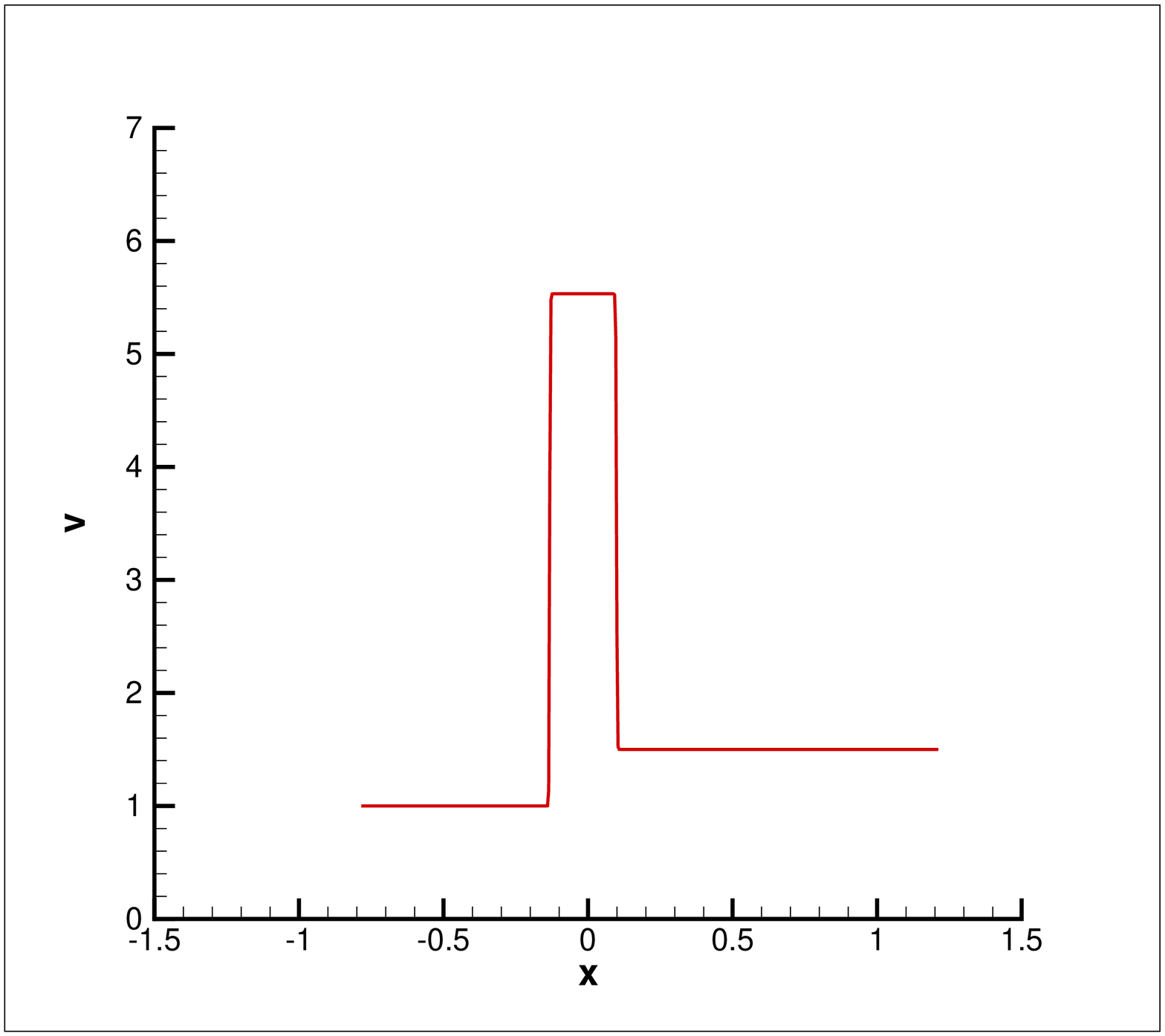}
\caption{Velocity and density for $\epsilon=0.3$ at $t=0.4$}
\end{center}
\end{figure}

\begin{figure}[th]
\begin{center}
\includegraphics[width=2.5in]{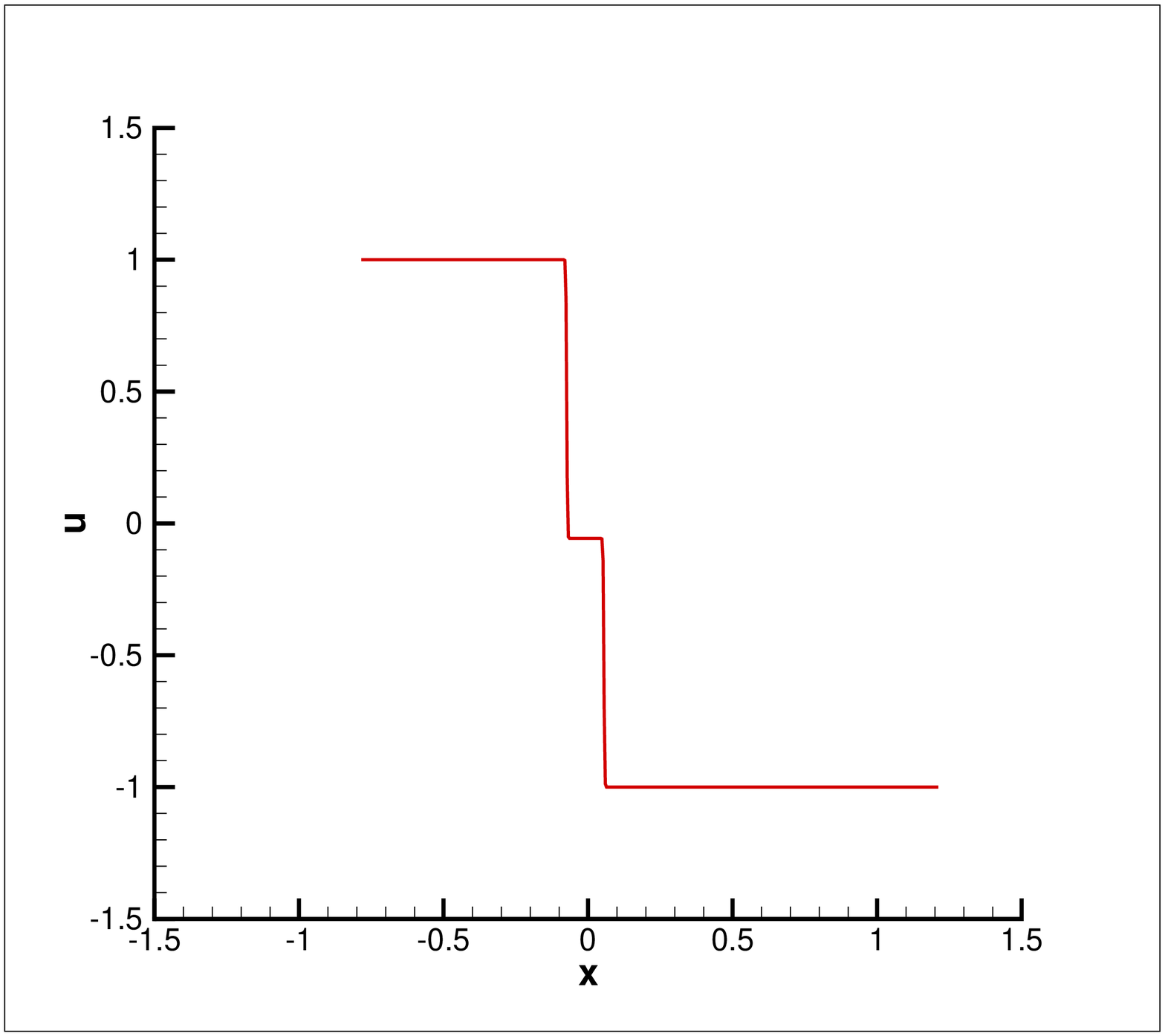}\hspace{0.50cm}\includegraphics[width=2.5in]{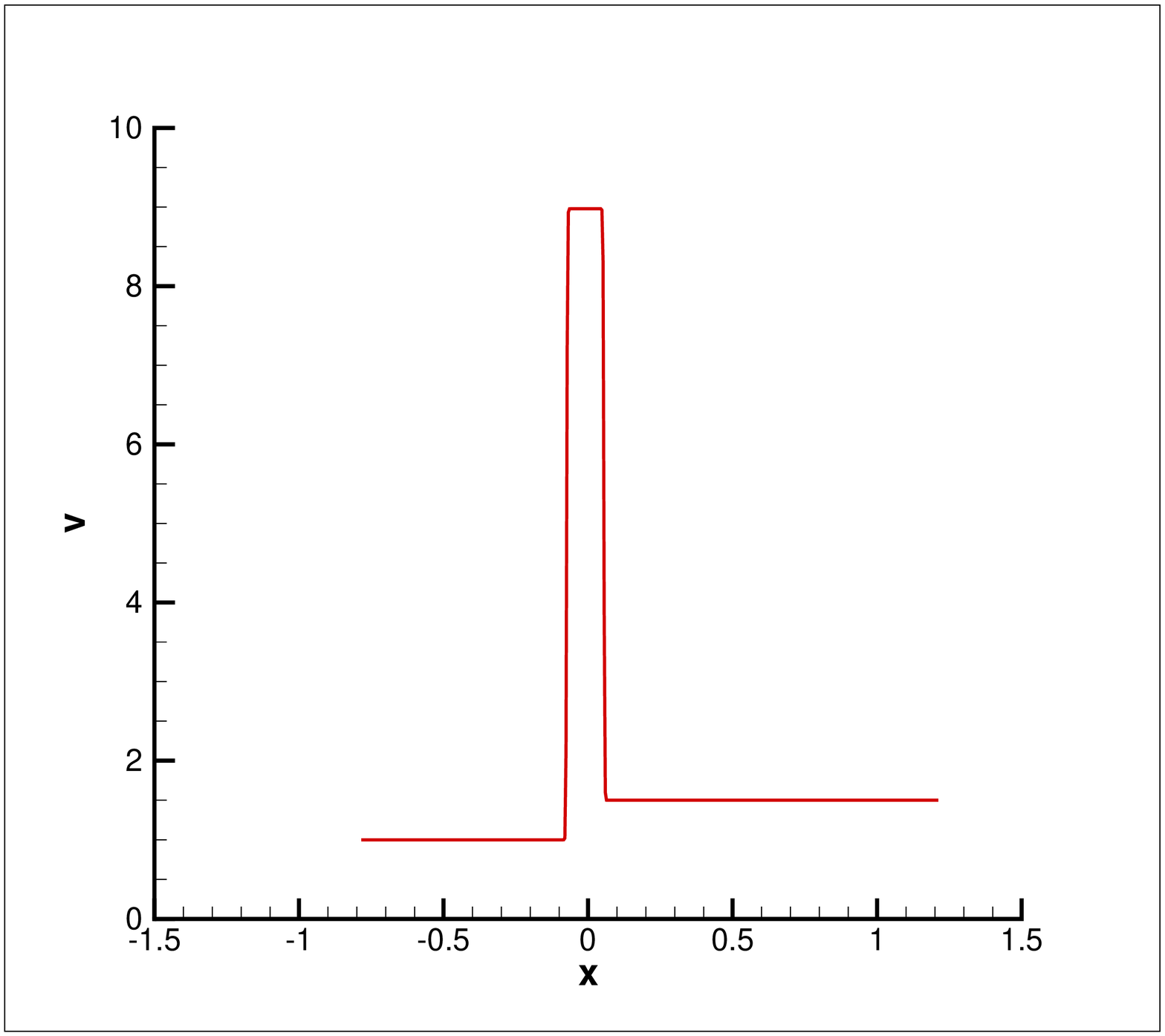}
\caption{Velocity and density for $\epsilon=0.15$ at $t=0.4$}
\end{center}
\end{figure}

\begin{figure}[th]
\begin{center}
\includegraphics[width=2.5in]{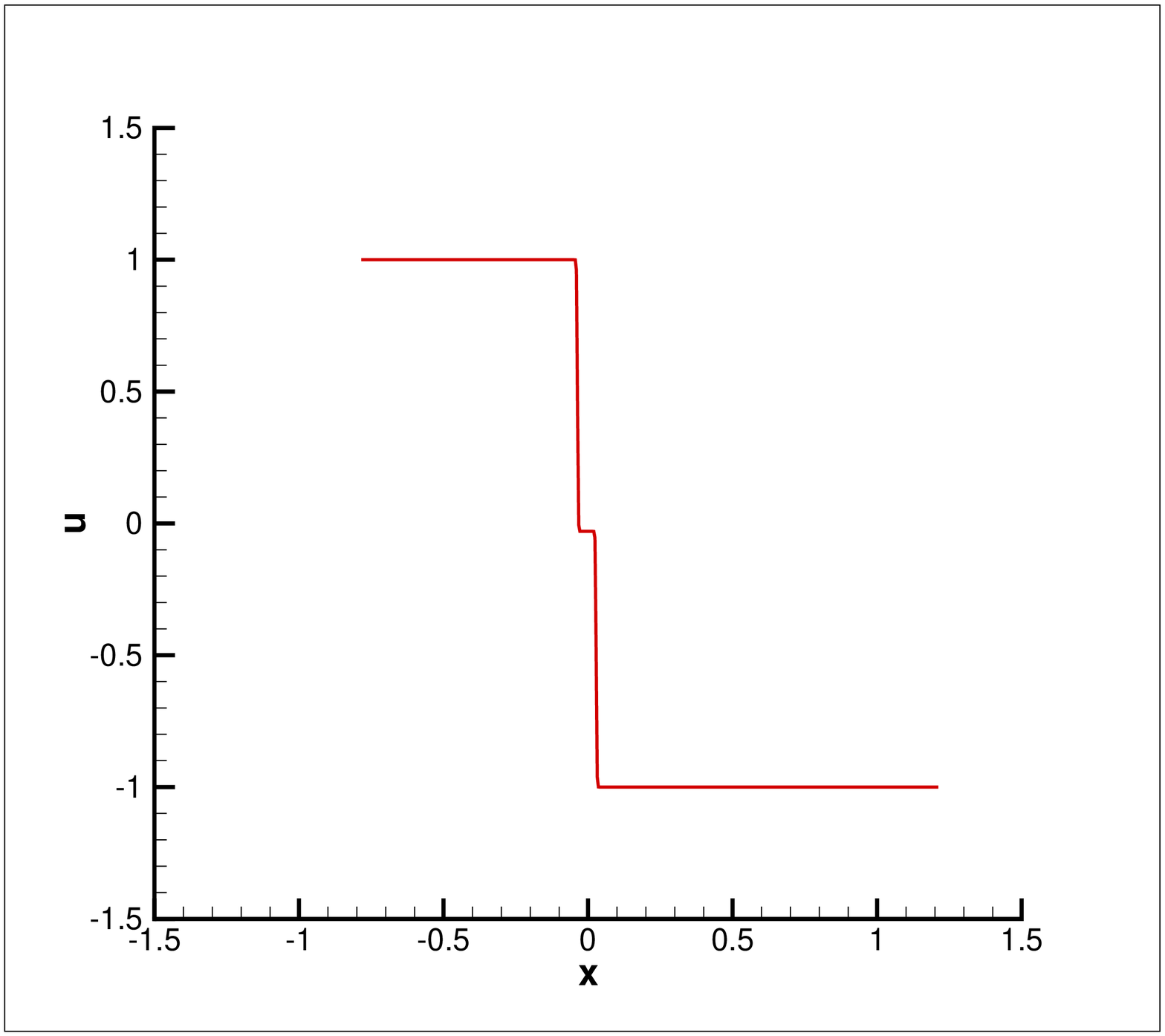}\hspace{0.50cm}\includegraphics[width=2.5in]{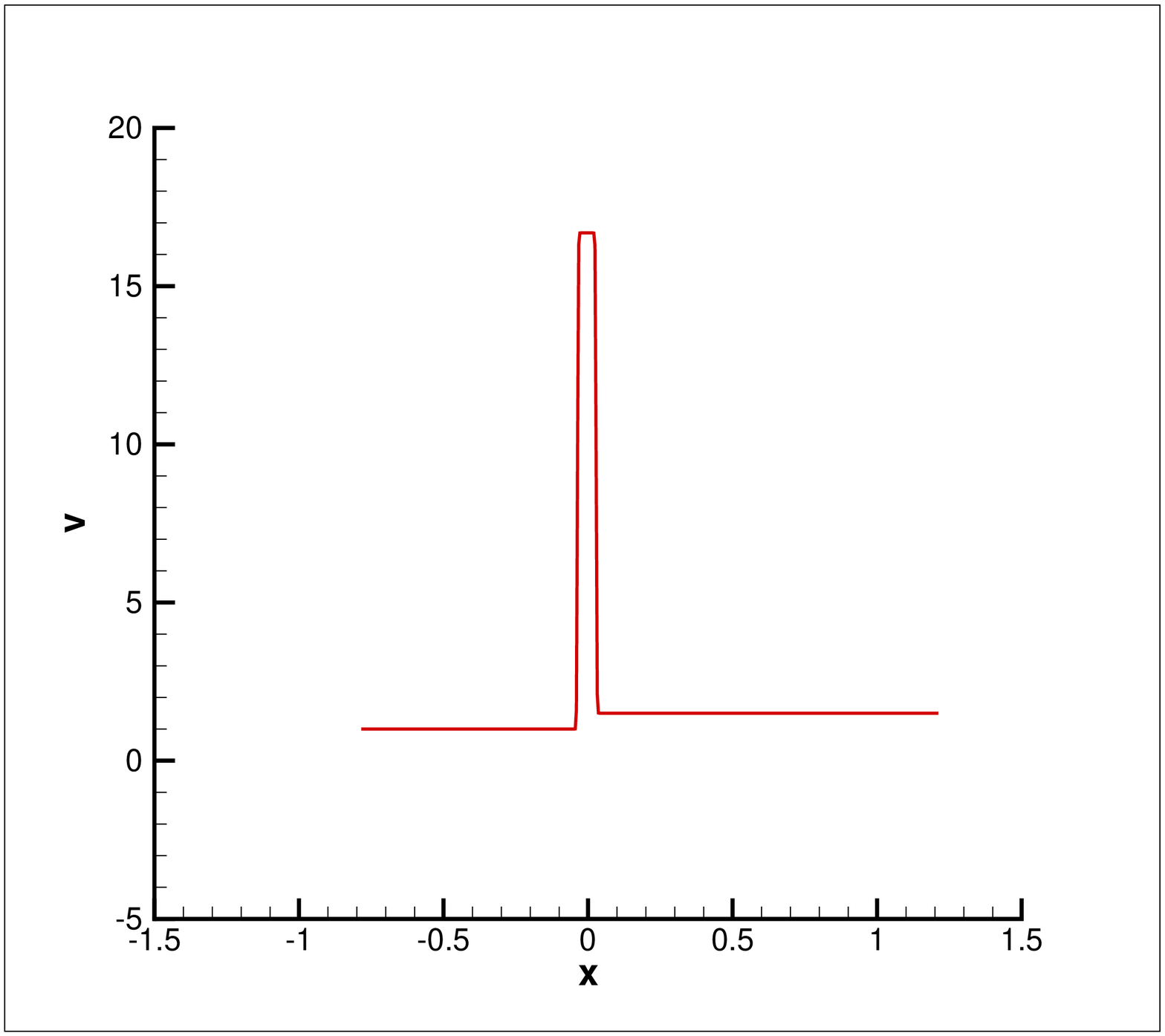}
\caption{Velocity and density for $\epsilon=0.07$ at $t=0.4$}
\end{center}
\end{figure}

\begin{figure}[th]
\begin{center}
\includegraphics[width=2.5in]{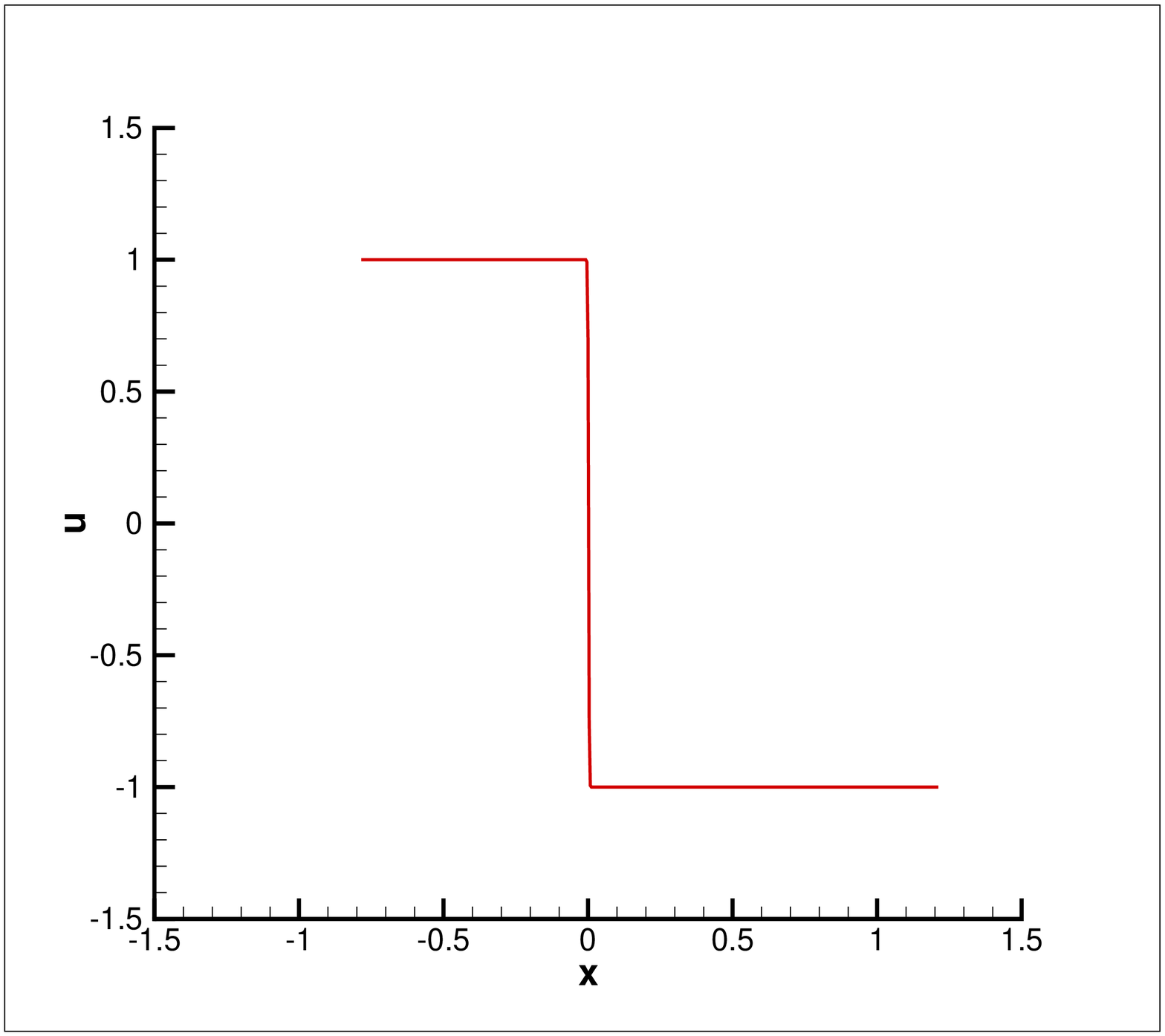}\hspace{0.50cm}\includegraphics[width=2.5in]{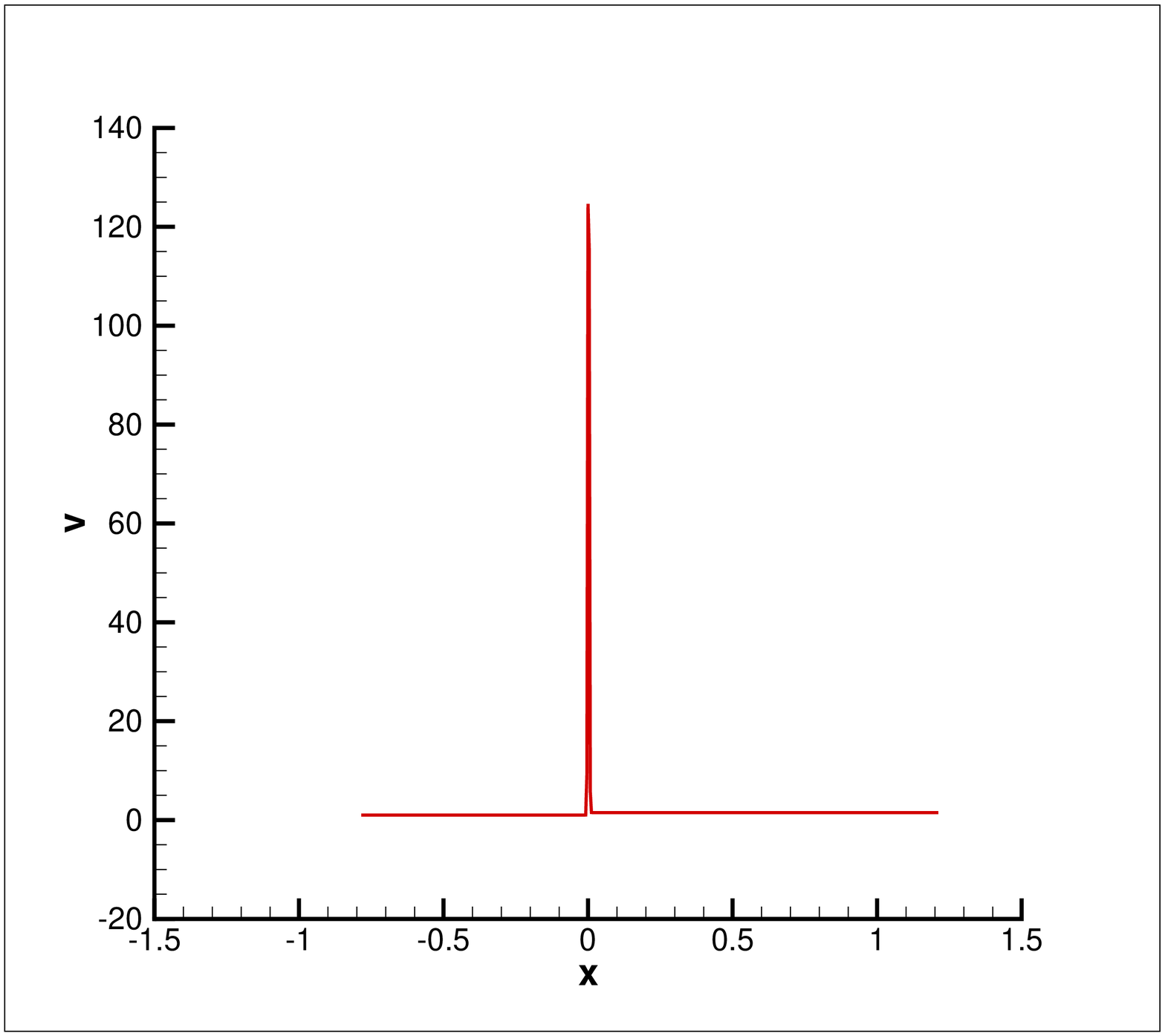}
\caption{Velocity and density for $\epsilon=0.001$ at $t=0.4$}
\end{center}
\end{figure}

One can observe clearly  from these above numerical results that,
when $\epsilon$ decreases, the location of the two shocks becomes
closer and closer, and the density of the intermediate state
increases dramatically, while the velocity is closer to a step
function. The numerical simulations are in complete agreement with
the theoretical analysis.


\begin{thebibliography}{00}

\bibitem{B}F. Bouchut, On zero-pressure gas dynamics, Advances in kinetic theory and computing, Series on
Advances in Mathematics for Applied Sciences, Vol. 22, World
Scientific, River Edge, NJ, 1994, pp. 171-190.

\bibitem{C-L-1} G. Chen, H. Liu, Formation of delta-shocks and vacuum states in the vanishing pressure limit of solutions to the
isentropic Euler equations,  SIAM Journal on Mathematical Analysis
34(2003)  925-938.

\bibitem{C-L-2} G. Chen, H. Liu, Concentration and cavitation in the vanishing pressure limit of solutions to the Euler equations for
nonisentropic fluids, Physica D Nonlinear Phenomena 189(2004)
141-165.


\bibitem{Cheng-Yang-1} H. Cheng, H. Yang,  Delta shock waves in chromatography
equations, Journal of Mathematical Analysis and Applications
380(2011)  475-485.

\bibitem{Cheng-1}H. Cheng, H. Yang, Riemann problem for the relativistic Chaplygin Euler
equations, Journal of Mathematical Analysis and Applications
381(2011) 17-26.



\bibitem{DS} V. Danilov,  V. Shelkovich, Dynamics of
propagation and interaction of $\delta$-shock waves in conservation
laws systems, Journal of Differential Equations  221(2005)  333-381.



\bibitem{D-M}  R. J. Diperna, A. Majda,  Reduced hausdorff dimension and
concentration-cancellation for two dimensional incompressible flow,
 Journal of the American Mathematical Society 1(1988)  59-59.



\bibitem{Fl} P. Le Floch, An existence and uniqueness result for two nonstrictly hyperbolic systems, in Nonlinear
Evolution Equations that Change Type, IMA 27 in Mathematics and its
Applications, Springer-Verlag, 1990


\bibitem{F-T} J. Fritz, B. Toth, Derivation of the Leroux system as the
hydrodynamic limit of a two-component lattice gas, Communications in
Mathematical Physics  249(2004) 1-27.





\bibitem{Gr} T. Gramchev, Entropy solutions to conservation laws with singular
initial data, Nonlinear Analysis, Theory, Methods \& Applications
24(1995) 721¨C733.




\bibitem{G-T} C. Greengard, E. Thomann, On diperna-majda concentration sets for
two-dimensional incompressible flow, Communications on pure \&
applied mathematics 41(1988)  295-303.





\bibitem{Hu}  J. Hu, A limiting viscosity approach to Riemann solutions containing
delta-shock waves for nonstrictly hyperbolic conservation laws,
Quarterly of Applied Mathematics  55(1997) 361-373.

\bibitem{J-S} P. Ji, C. Shen, Construction of the global solutions to the perturbed
Riemann problem for the Leroux system, Advances in Mathematical
Physics, Volume 2016, Article ID 4808610, 13 pages.


\bibitem{Jiang} G. Jiang, E. Tadmor,  Non-oscillatory central schemes for multidimensional hyperboloic conservation laws,
 SIAM Journal on Scientific Computing 19(1998)  1892-1917.




\bibitem{K} K. T. Joseph, A Riemann problem whose viscosity solutions contain delta-measures, Asymptotic
Analysis 7, 105-120 (1993).



\bibitem{KK}B. L. Keyfitz and H. C. Kranzer, A viscosity approximation to system of conservation laws with
no classical Riemann solution in Nonlinear Hyperbolic Problems,
Lecture Notes in Mathematics, Vol. 1042, Springer-Verlag, Berlin/New
York, 1989.



\bibitem{Ko}  D. J. Korchinski,  Solutions  of a Riemann  problem for a 2x2
system  of conservation  laws pos- sessing  classical  solutions,
Adelphi  University  Thesis,  1977.


\bibitem{Le}  A. Y. Leroux, Approximation des systems hyperboliques, in
Cours et Seminaires INRIA, Problems Hyperboliques, Rocquen-court,
1981.



\bibitem{LJQ} J. Li, Note on the compressible Euler equations with
zero temperature, Applied Mathematics Letters 14(2001) 519-523.


\bibitem{L-M}  Y. G. Lu, I. Mantilla, and L. Rendon, Convergence of
approximated solutions to a nonstrictly hyperbolic system, Advanced
Nonlinear Studies  1(2001) 65-79.

\bibitem{Lu}  Y. G. Lu, Global entropy solutions of Cauchy problem for the
Le Roux system, Applied Mathematics Letters 60(2016) 61-66.





\bibitem{D-MN} Darko Mitrovic, Marko Nedeljkov   Delta shock waves as a limit of shock
waves, Journal of Hyperbolic Differential Equations 4 (2007) 1-25.


\bibitem{N-T} H. Nessyahu, E. Tadmor, Non-oscillatory central differencing for hyperbolic conservation laws,   Journal of
Computational Physics 87(1990)  408-463.






\bibitem{Se-1} D. Serre,  Systems of Conservation Laws, Vol. 1-2, Cambridge: Cambridge University Press,
2000.



\bibitem{Se-2} D. Serre,   Existence globale de solutions faibles sous une
hypothse unilaterale pour un systme hyperbolique non linaire,
Quarterly of Applied Mathematics 46(1988) 157-167.




\bibitem{VSh}V. Shelkovich,  The Riemann problem admitting $\delta$, $\delta'$-shocks,
and vacuum states (the vanishing viscosity approach), Journal of
Differential Equations 231(2006) 459-500.


\bibitem{Sh}  C. Shen, M. Sun, Formation of delta shocks and vacuum states in the
vanishing pressure limit of Riemann solutions to the perturbed
Aw-Rascle model,  Journal of Differential Equations  249(2010)
3024-3051.


\bibitem{SZh} W. Sheng,  T. Zhang, The Riemann problem for the transportation equations in gas dynamics,
Mem. Amer. Math. Soc. 137 (1999).






\bibitem{Ta} D. Tan, Riemann problem for hyperbolic systems of conservation laws
with no classical wave solutions, Quarterly of Applied Mathematics
51(1993) 756-776.



\bibitem{T-Z-Z} D.  Tan,  T. Zhang, Y.  Zheng,  Delta shock waves as
limits of vanishing viscosity for hyperbolic systems of conservation
laws, Journal of Differential Equations  112(1994) 1-32.





\bibitem{Tan-Zhang} D.  Tan,  T. Zhang,  Two-dimensional Riemann problem for
a hyperbolic system of nonlinear conservation laws I. Four-J cases,
II. Initial data involving some rarefaction waves, Journal of
Differential Equations 111(1994)  203-282.





\bibitem{Te} B. Temple, Systems of conservation laws with invariant
submanifolds, Transactions of the American Mathematical Society vol.
280(1983) 781-795.




\bibitem{F-V} B. Toth, B. Valko, Perturbation of singular equilibria of hyperbolic two-component
systems: A universal hydrodynamic limit, Communications in
Mathematical Physics 256(2005) 111-157.



\bibitem{Yang-1}H. Yang, Riemann Problems for a Class of Coupled Hyperbolic Systems
of Conservation Laws,  Journal of Differential Equations   159(1999)
447-484.


\bibitem{Yang-Liu-1} H. Yang,  J. Liu,  Delta-shocks and vacuums in pressureless gas
dynamics by the flux approximation, Science China Mathematics
58(2015)  2329-2346.



\bibitem{Yang-Liu-2} H. Yang,  J. Liu,   Concentration and cavitation in the Euler equations for
 nonisentropic fluids with the flux approximation,  Nonlinear Analysis   123-124(2015) 158-177.

\bibitem{YSh-1} G. Yin, W. Sheng, Delta shocks and vacuum states in vanishing
pressure limits of solutions to the relativistic Euler equations for
polytropic gases, Journal of Mathematical Analysis and Applications
355 (2009) 594-605.


\bibitem{YSh-2} G. Yin, W. Sheng, Delta wave formation and vacuum state in vanishing
pressure limit for system of conservation laws to relativistic fluid
dynamics, Zamm Journal of Applied Mathematics \& Mechanics 95(2013)
49-65.










\end{thebibliography}
\end{document}